\documentclass[12pt,reqno]{article}
\usepackage
[colorlinks=true, linkcolor=webgreen, filecolor=webbrown,
citecolor=webgreen]{hyperref}

\usepackage{amsmath}
\usepackage{amssymb}
\usepackage{color}
\usepackage{epsf}

\begin{document}
\makeatletter

\begin{center}
\epsfxsize=10in
\end{center}

\def\endofproofmark{$\Box$}

\def\endofproofmark{$\Box$}
\begin{center}
\vskip 1cm {\LARGE\bf An Extension of Alzer's Inequality} \vskip
8mm {\LARGE\bf by Convexity} \vskip 1cm \large J. Rooin
\\
\vskip .5cm
Department of Mathematics\\
Institute for Advanced
Studies in Basic Sciences\\
P.O. Box 45195-1159\\
Zanjan, Iran\\
\href{mailto:rooin@iasbs.ac.ir}{\tt rooin@iasbs.ac.ir}\\
\end{center}



\date{}
\newtheorem{theo}{Theorem}[section]
\newtheorem{prop}[theo]{Proposition}
\newtheorem{rem}{Remark}
\newtheorem{lemma}[theo]{Lemma}
\newtheorem{cor}[theo]{Corollary}
\newtheorem{problem}{Problem}
\def\frameqed{\framebox(5.2,6.2){}}
\def\deshqed{\dashbox{2.71}(3.5,9.00){}}
\def\ruleqed{\rule{5.25\unitlength}{9.75\unitlength}}
\def\myqed{\rule{8.00\unitlength}{12.00\unitlength}}
\def\qed{\hbox{\hskip 6pt\vrule width 7pt height11pt depth1pt\hskip 3pt}
\bigskip}
\newenvironment{proof}{\trivlist\item[\hskip\labelsep{\bf Proof}.]}{\hfill
 $\frameqed$ \endtrivlist}
\newcommand{\COM}[2]{{#1\choose#2}}

\thispagestyle{empty} \null \addtolength{\textheight}{1cm}

\begin{abstract}
In this article, we obtain two interesting general inequalities
concerning Riemman sums of convex functions, which in particular,
sharpen Alzer's inequality and give a suitable converse for it.
\end{abstract}

\bigskip
\hrule
\bigskip



\bigskip
\hrule
\bigskip

\section{Introduction}
In [1], H. Alzer proved the following inequality
\begin{eqnarray}
\frac{n}{n+1}\le\left((n+1)\sum_{i=1}^ni^r/n\sum_{i=1}^{n+1}i^r
\right)^{1/r},
\end{eqnarray}
where $r$ is a positive real and $n$ is a natural number. In other
words, the Riemman sums
$\frac{1}{n}\sum_{i=1}^n\left(\frac{i}{n}\right)^r~(n=1,2,\cdots)$
of the function $x^r$ is a decreasing sequence. The proof of Alzer
[1] is technical, but quite complicated. So, in several articles
Alzer's proof has been simplified, and also in many others, this
inequality has been extended; see e.g. [2-4]. \\
In this article, using some trivial facts about convex functions,
we obtain some valuable results concerning special kinds of
Riemman sums of convex functions, from which Alzer's inequality
with an its converse are handled at once.
\section{Main Results}
Throughout this section, we suppose that
$f:[a,b]\rightarrow\mathbb{R}$ is an arbitrary function on a
closed interval $[a,b]$, and put
\begin{eqnarray}
A_n=\frac{b-a}{n}\sum_{i=1}^n
f\left(x_i^{(n)}\right)\hspace{1cm}\mbox{and}\hspace{1cm}
B_n=\frac{b-a}{n}\sum_{i=0}^{n-1} f\left(x_i^{(n)}\right),
\end{eqnarray}
where
$$
x_i^{(n)}=a+i\frac{b-a}{n}\hspace{1.5cm}(i=0,1,\cdots,n;~n=1,2,\cdots)
$$
When emphasizing, we write $A_n(f)$ instead of $A_n$, and so on.
In the following theorem, we obtain some recursive inequalities
concerning $A_n$ and $B_n$, which as a corollary, give Alzer's
inequality and a converse for it.\\\\
{\bf Theorem 2.1} With the above assumptions, if $f$ is convex,
then we have
\begin{eqnarray} A_{n+1}+\frac{1}{n(n+2)}\left[A_{n+1}-(b-a)f(a))\right]\le
A_n\le A_{n+1}+\frac{1}{n^2}\left[(b-a)f(b)-A_{n+1}\right]
\end{eqnarray}
and
\begin{eqnarray} B_{n+1}+\frac{1}{n(n+2)}\left[B_{n+1}-(b-a)f(b))\right]\le
B_n\le B_{n+1}+\frac{1}{n^2}\left[(b-a)f(a)-B_{n+1}\right].
\end{eqnarray}
Moreover,
\begin{eqnarray}
A_{n+1}\le
(b-a)\left[\frac{n}{2(n+1)}f(a)+\frac{n+2}{2(n+1)}f(b)\right]
\end{eqnarray}
and
\begin{eqnarray}
B_{n+1}\le
(b-a)\left[\frac{n+2}{2(n+1)}f(a)+\frac{n}{2(n+1)}f(b)\right].
\end{eqnarray}
If $f$ is concave, all the above inequalities reverse.
Moreover, all these inequalities are strict in the case of strict convexity or concavity. \\
{\it
Proof}. Since
$x_i^{(n+1)}=\frac{i}{n+1}x_{i-1}^{(n)}+\frac{n+1-i}{n+1}x_i^{(n)}~(1\le
i\le n)$, by Jensen's inequality, we have
$$
f(x_i^{(n+1)})\leq\frac{i}{n+1}f(x_{i-1}^{(n)})+\frac{n+1-i}{n+1}f(x_i^{(n)})\hspace{1.5cm}(1\le
i\le n),
$$
which by summing them up from $i=1$ to $i=n$, with some
calculations, we get the left
hand side of (3).\\
Similarly, since
$x_{i-1}^{(n)}=\frac{n+1-i}{n}x_{i-1}^{(n+1)}+\frac{i-1}{n}x_i^{(n+1)}~(1\le
i\le n+1)$, we have
$$
f(x_{i-1}^{(n)})\leq
\frac{n+1-i}{n}f(x_{i-1}^{(n+1)})+\frac{i-1}{n}f(x_i^{(n+1)})\hspace{1.5cm}(1\le
i\le n+1),
$$
which by summing them up from $i=1$ to $i=n+1$, with some
calculations, we get the right hand
side of (4).\\
The right hand side of (3) and the left hand side of (4) are
obtained from the right hand side of (4) and the left hand side of
(3) respectively, by considering
$$
A_k-B_k=\frac{(b-a)}{k}[f(b)-f(a)]\hspace{1.5cm}(k=1,2,\cdots).$$
The inequalities (5) and (6) are trivially obtained by comparing
the left and right hand sides
of (4) and (5) with each other respectively.\\\\
{\bf Corollary 2.2.} If $f$ is an increasing convex or concave
function on $[a,b]$, then
\begin{eqnarray}
A_{n+1}\leq A_n\hspace{1cm}\mbox{and}\hspace{1cm}B_n\leq
B_{n+1}\hspace{1cm}(n=1,2,\cdots).
\end{eqnarray}
{\it Proof}. Since $f$ is increasing, for each $k$ we have
\begin{eqnarray}
(b-a)f(a)\leq A_k\leq(b-a)f(b)
\end{eqnarray}
and
\begin{eqnarray}
(b-a)f(a)\leq B_k\leq(b-a)f(b).
\end{eqnarray}
So, if $f$ is convex, using (8), (9), the left hand side
of (3) and the right hand side of (4), we get (7).\\
Now, if $f$ is concave, then $-f$ is convex, and (7) follows from
the right hand side of (3) and the left hand side of (4) applying
to $-f$, by using (8) and (9) and taking into consideration that
$$
A_k(-f)=-A_k(f)\hspace{1cm}\mbox{and}\hspace{1cm}B_k(-f)=-B_k(f)\hspace{1cm}(k=1,2,\cdots).
$$
\\
{\bf Corollary 2.3.} If $r\geq 1$, then
\begin{eqnarray}
\frac{n}{n+1}\left(1+\frac{1}{n(n+2)}\right)^{1/r}\leq
\left((n+1)\sum_{i=1}^ni^r/n\sum_{i=1}^{n+1}i^r \right)^{1/r}
\end{eqnarray}
$$\leq\frac{n}{n+1}\left(1+\frac {(n+1)^{r+1}-\sum_{i=1}^{n+1}i^r}
{n^2\sum_{i=1}^{n+1}i^r} \right)^{1/r}.
$$
If $0<r\leq 1$  the inequalities in (10) reverse. \\
Obviously, these give us a refinement and a reverse of Alzer's
inequality
(1).\\
{\it Proof}. If $r\geq 1$, the function $f(x)=x^r~(x\geq 0)$ is
convex and so (10) follows from (3), by taking $a=0$ and $b=1$.\\
If $0<r\leq 1$ the function $f$ is concave and so $-f$ is convex.
So, the inequalities in (3), and
therefore the inequalities in (10) reverse .\\
\begin{center}
{REFERENCES}
\end{center}
\begin{enumerate}
\item H. Alzer, On an inequality of H. Minc and L. Sathre, {\it J. Math. Anal. Appl.}, {\bf 179} (1993),
396-402.
\item F. Qi, Generalization of
H. Alzer's inequality, {\it J. Math. Anal. Appl.} {\bf 240}
(1999), no. 1, 294--297.
\item Elezovi\'c, N.; Pe\v cari\'c, J. On Alzer's
inequality. {\it J. Math. Anal. Appl.} {\bf 223} (1998), no. 1,
366--369.
\item J. Sandor, On an inequality of Alzer, {\it J. Math. Anal. Appl.} {\bf 192} (1995), no. 3,
1034-1035.
\end{enumerate}
\end{document}